\documentclass[a4paper,leqno,11pt,final]{amsproc}
\usepackage{amssymb,latexsym,amsxtra,amscd,amsfonts,amsthm,amsmath,verbatim}

\numberwithin{equation}{section}

\begin{document}
 \theoremstyle{plain}
   \newtheorem{theorem}                    {Theorem}       [section]
   \newtheorem{lemma}      [theorem]       {Lemma}
   \newtheorem{corollary}  [theorem]       {Corollary}
   \newtheorem{proposition}[theorem]       {Proposition}
 \theoremstyle{definition}
   \newtheorem{definition} [theorem]       {Definition}
   \newtheorem{conjecture} [theorem]       {Conjecture}
 \theoremstyle{remark}
   \newtheorem{remark}     [theorem]       {Remark}

\catcode`@=11
\atdef@ I#1I#2I{\CD@check{I..I..I}{\llap{$\m@th\vcenter{\hbox
  {$\scriptstyle#1$}}$}
  \rlap{$\m@th\vcenter{\hbox{$\scriptstyle#2$}}$}&&}}
\atdef@ E#1E#2E{\ampersand@
  \ifCD@ \global\bigaw@\minCDarrowwidth \else \global\bigaw@\minaw@ \fi
  \setboxz@h{$\m@th\scriptstyle\;\;{#1}\;$}%
  \ifdim\wdz@>\bigaw@ \global\bigaw@\wdz@ \fi
  \@ifnotempty{#2}{\setbox@ne\hbox{$\m@th\scriptstyle\;\;{#2}\;$}%
    \ifdim\wd@ne>\bigaw@ \global\bigaw@\wd@ne \fi}%
  \ifCD@\enskip\fi
    \mathrel{\mathop{\hbox to\bigaw@{}}%
      \limits^{#1}\@ifnotempty{#2}{_{#2}}}%
  \ifCD@\enskip\fi \ampersand@}
\catcode`@=\active

\newcommand{\Aut}{\operatorname{Aut}}
\newcommand{\Hom}{\operatorname{Hom}}
\newcommand{\End}{\operatorname{End}}
\newcommand{\Ext}{\operatorname{Ext}}
\newcommand{\Gal}{\operatorname{Gal}}
\newcommand{\Pic}{\operatorname{Pic}}
\newcommand{\ord}{\operatorname{ord}}
\newcommand{\Spec}{\operatorname{Spec}}
\newcommand{\trdeg}{\operatorname{trdeg}}
\newcommand{\holim}{\operatornamewithlimits{holim}}
\newcommand{\im}{\operatorname{im}}
\newcommand{\coim}{\operatorname{coim}}
\newcommand{\coker}{\operatorname{coker}}
\newcommand{\gr}{\operatorname{gr}}
\newcommand{\id}{\operatorname{id}}
\newcommand{\Br}{\operatorname{Br}}
\newcommand{\cd}{\operatorname{cd}}
\newcommand{\tr}{\operatorname{tr}}
\newcommand{\CH}{\operatorname{CH}}
\renewcommand{\lim}{\operatornamewithlimits{lim}}
\newcommand{\ilim}{\operatornamewithlimits{lim}}
\newcommand{\colim}{\operatornamewithlimits{colim}}
\newcommand{\rk}{\operatorname{rank}}
\newcommand{\codim}{\operatorname{codim}}
\newcommand{\NS}{\operatorname{NS}}
\newcommand{\TC}{\operatorname{TC}}
\newcommand{\TR}{\operatorname{TR}}
\renewcommand{\TH}{\operatorname{TH}}
\newcommand{\cone}{{\rm cone}}
\newcommand{\f}{{\mathcal F}}
\newcommand{\p}{{\mathcal P}}
\newcommand{\N}{{\mathbb N}}
\newcommand{\A}{{\mathbb A}}
\newcommand{\Z}{{{\mathbb Z}}}
\newcommand{\Q}{{{\mathbb Q}}}
\newcommand{\R}{{{\mathbb R}}}
\renewcommand{\H}{{{\mathbb H}}}
\renewcommand{\P}{{{\mathbb P}}}
\newcommand{\K}{{K}}
\newcommand{\Zp}{{\Z_p}}
\newcommand{\Qp}{{\Q_p}}
\newcommand{\Zpn}{{\Z/p^n}}
\newcommand{\F}{{{\mathbb F}}}
\newcommand{\m}{{\mathfrak m}}
\newcommand{\Sm}{\text{\rm Sm}}
\newcommand{\Sch}{\text{\rm Sch}}
\newcommand{\et}{{\text{\rm\'et}}}
\newcommand{\Zar}{{\text{\rm Zar}}}
\newcommand{\Nis}{{\text{\rm Nis}}}
\newcommand{\PreShv}{\text{\rm PreShv}}
\newcommand{\Div}{\operatorname{Div}}
\renewcommand{\div}{\operatorname{div}}
\newcommand{\corank}{\operatorname{corank}}
\renewcommand{\O}{{\mathcal O}}
\renewcommand{\p}{{\mathfrak p}}

\title{The cyclotomic trace map and values of zeta functions}
\author{Thomas Geisser}
\thanks{Supported in part by NSF grant No. 0300133, 
and the Alfred P. Sloan Foundation}
\address{University of Southern California, 
Dep. of Math., KAP 108, 3620 Vermont Av.,
Los Angeles 90089}
\email{geisser\char`\@math.usc.edu}
\begin{abstract}
We show that the cyclotomic trace map for smooth varieties
over number rings can be interpreted as a regulator map 
and hence are related to special values of $\zeta$-functions.
\end{abstract}
\maketitle

\section{Introduction}
The purpose of this paper is to show that the results of
\cite{ichlarsmixed} can be used to relate the cyclotomic trace map
from \'etale $K$-theory to topological cyclic homology
\begin{equation*}
\tr_i:K_i^\et(X,\Z_p)\to \TC_i(X;p,\Z_p)
\end{equation*}
to arithmetic invariants if $X$ is regular scheme, flat and proper over a 
number ring $\O_S$, and with good reduction at $p$. 
The main result of \cite{ichlarsmixed} implies that the map 
$\tr_i$ can be identified with the localization map
\begin{equation}\label{map}
\tr_i: K_i^\et(X,\Z_p)\to K_i^\et(X\times_\Z\Z_p,\Z_p).
\end{equation}
Both sides of \eqref{map} admit a hypercohomology spectral sequence of the 
form 
$$E_2^{s,t}=H^s_{cont}(X,({\mathcal K}/p^\cdot)_{-t}) \Rightarrow 
K_{-s-t}^\et(X,\Z_p),$$
but the $E_2$-term is hard to control because
the \'etale $K$-theory sheaf $({\mathcal K}/p^r)_i$ is not known.
In order to overcome this problem,
we compare the map \eqref{map} to the map 
\begin{equation} \label{map1}\textstyle
K_i^\et(X\times_\Z\Z[\frac{1}{p}],\Z_p)
\to K_i^\et(X\times_\Z\Q_p,\Z_p).
\end{equation}
Using either the Lichtenbaum-Quillen conjecture or a result
of Thomason, one can identify the maps \eqref{map} and \eqref{map1}
if one assumes that $i>d=\dim X$ or $i\geq \frac{8}{3}(d+2)(d+3)(d+4)-14$,
respectively. Since $p$ is invertible in \eqref{map1}, the \'etale
$K$-theory sheaf can be identified, and the map \eqref{map1}
is the map on the abutments of the spectral sequences
\begin{align*}
\begin{CD}
\textstyle
E_2^{s,t}=H^s(X\times_\Z\Z[\frac{1}{p}],\Z_p(-\frac{t}{2})) 
@E\Rightarrow EE   \textstyle
K_{-s-t}^\et(X\times_\Z\Z[\frac{1}{p}],\Z_p)\\ 
@VVV @VVV \\
\textstyle E_2^{s,t}=H^s(X\times_\Z\Q_p,\Z_p(-\frac{t}{2})) 
@E\Rightarrow EE  
\textstyle K_{-s-t}^\et(X\times_\Z\Q_p,\Z_p).
\end{CD}
\end{align*}
Thus we can relate the map \eqref{map1} to maps between \'etale 
cohomology groups, which in turn are related to special values
of $L$-functions.

In the second half of the paper we give concrete calculations in case
$X$ is the ring of integers $\Spec \O$ of a number field $K$. 
For $j=1,2$, the trace map $tr_{2i-j}$ can be identified with the map 
$$H^j(G_\Sigma,\Z_p(i)) \to \prod_{\p|p} H^j(K_\p,\Z_p(i)),$$
where $G_\Sigma$ is Galois group of the maximal extension of $K$
which is unramified outside of $p$ and infinity.
This map has been study in Iwasawa theory, and we translate results
of Iwasawa theory into statements about the trace map.
For example, we show that if $K/\Q$ is a totally real field, 
unramified at $p$, and $i\not\equiv 1\mod p-1$ is an odd integer, 
then the trace map 
$$\tr_{2i-1}'  :K_{2i-1}(\O)\otimes \Z_p/\text{tors} \to 
\TC_{2i-1}(\O,\Z_p) .$$
is a map between free $\Z_p$-modules of rank $d=[K:\Q]$.
Its cokernel is finite if and only if a conjecture of Schneider
holds, and in this case the order of the cokernel is related
to the $p$-adic $L$-function as follows:
$$ |H^2(G_\Sigma,\Z_p(i))|\cdot|\coker \tr_{2i-1}'|=
|L_p(K,\omega^{1-i},i)|_p^{-1}.$$

{\it Convention:} All cohomology groups are \'etale cohomology 
in case of schemes, and Galois cohomology groups in case of
fields.

\section{Preliminaries} 
We recall some facts on algebraic $K$-theory and topological
cyclic homology, see \cite{handbook}, \cite{ichlars} and \cite{larsicm}.

\subsection{K-theory }
For every henselian pair $(A,I)$ such that $m$ is invertible in $A$,
and for all $i\geq 0$, we have the isomorphism of Gabber \cite{gabberrig} 
and Suslin \cite{susdvr}
\begin{equation*}
K\sb i(A, {\mathbb Z}/m) 
\stackrel{\sim}{\longrightarrow}K\sb i(A/I,{\mathbb Z}/m).
\end{equation*}
Together with the calculation of the $K$-theory of an
algebraically closed field \cite{susdvr} by Suslin, this implies that
on every scheme $X$ such that $m$ is invertible on $X$, the
$K$-theory sheaf with coefficients for the \'etale topology
can be identified as follows
\begin{equation}\label{modmsheaf}
({\mathcal K}/m)_n=
\begin{cases}
\mu_m^{\otimes \frac{n}{2}} &n\geq 0 \;\text{even},\\
0 & n \;\text{odd}.
\end{cases}
\end{equation}

Let
$R$ be a local ring, such that $(R,pR)$ is a henselian pair,
and such that $p$ is not a zero divisor. Then \cite{ichlarsmixed} 
the reduction map
\begin{equation}
\label{contin}
K_i(R,\Z/p^r)\to \{K_i(R/p^s,\Z/p^r)\}_s
\end{equation}
is an isomorphism of pro-abelian groups. This generalizes
the result of Suslin and Panin \cite{susdvr,panin} for
$R$ a henselian valuation ring  of mixed characteristic $(0,p)$.

For a presheaf of spectra $\f$ on a site $X_\tau$, and a covering
${\mathcal U}=\{U_i\}$ of $X$, Thomason \cite[Def. 1.9, 1.33]{thomasonetale} 
defines the 
\v Cech hypercohomology spectrum $\check{\mathbb H}^\cdot({\mathcal U},\f)$ 
and the sheaf hypercohomology spectrum ${\mathbb H}^\cdot(X_\tau,\f) $.
There are natural augmentation maps 
$\tau:\f(X)\to \check {\mathbb H}^\cdot({\mathcal U},\f)$ and
$\eta:\f(X)\to {\mathbb H}^\cdot(X_\tau,\f)$.
If $\tau$ is the Zariski or Nisnevich topology on a noetherian
scheme $X$ of finite Krull dimension, then it is a theorem
of Brown-Gersten \cite{browngersten} and Nisnvevich \cite{nisnevich}, 
respectively, that the augmentation map
$\eta:K(X)\to {\mathbb H}^\cdot(X_\tau,K)$ is a homotopy equivalence.
Moreover, the \v Cech hypercohomology of the sheaf
hypercohomology agrees with the sheaf hypercohomololgy
\cite[Cor. 1.47]{thomasonetale} 
\begin{equation}\label{cechissheaf}
{\mathbb H}^\cdot(X_\Zar,\f)\cong
\check {\mathbb H}^\cdot({\mathcal U},{\mathbb H}^\cdot(-,\f)).
\end{equation}
One important feature of ${\mathbb H}^\cdot(X_\tau,\f)$ is that it comes
equipped with a spectral sequence \cite[Prop. 1.36]{thomasonetale}
\begin{equation}
\label{thomss} E_2^{s,t}=H^s(X_\tau,\tilde\pi_{-t}\f) \Rightarrow
\pi_{-s-t}{\mathbb H}^\cdot(X_\tau,\f),
\end{equation}
where
$\tilde\pi_i\f$ is the sheaf associated to the presheaf of homotopy groups
$U\mapsto \pi_i\f(U)$. If $X_\tau$ has finite cohomological dimension,
then the spectral sequence converges.

For a pro-presheaf $\f^\cdot$ 
of spectra on $X_\tau$, one defines the hypercohomology spectrum
${\mathbb H}^\cdot(X_\tau,\f^\cdot):=\holim_r{\mathbb H}^\cdot(X_\tau,\f^r)$.
For a (complex of) sheaves of abelian groups $A^\cdot$ on $X_\tau$,
Jannsen \cite{jan} defines continuous cohomology groups 
$H^j_{cont}(X_\tau,A^\cdot)$ as the derived functors of the functor
$ A^\cdot \mapsto \ilim_r\Gamma(X,A^r)$, and we get a spectral sequence 
\cite{ichlars}
\begin{equation}
\label{thomsscont} E_2^{s,t}=H^s_{cont}(X_\tau,\tilde\pi_{-t}\f^\cdot)
\Rightarrow \pi_{-s-t}{\mathbb H}^\cdot(X_\tau,\f^\cdot).
\end{equation}

If $X_\et$ is the small \'etale site
of the scheme $X$, then we write $K^\et_i(X,\Z_p)$ for the
homotopy groups
$\pi_i\holim_r {\mathbb H}^\cdot(X_\et,K/p^r)$.
If $p$ is invertible on $X$, we write 
$H^i(X,\Z_p(n))$ for $H^i_{cont}(X_\et,\mu_{p^\cdot}^{\otimes n})$.
In view of \eqref{modmsheaf} the spectral
sequence \eqref{thomsscont} takes the form
\begin{equation}\label{qwer} \textstyle
E_2^{s,t}=H^s(X,\Z_p(-\frac{t}{2}))\Rightarrow
K_{-s-t}^\et(X,\Z_p).
\end{equation}

The Lichtenbaum-Quillen conjecture states that for $i$
greater than the cohomological dimension of $X$,
the canonical map from $K$-theory to \'etale $K$-theory
$$K_i(X,\Z_p)\to K_i^\et(X,\Z_p)$$
is an isomorphism. The Lichtenbaum-Quillen conjecture
is a consequence of the Beilinson-Lichtenbaum conjecture,
whose proof has been announced by Voevodsky \cite{voebk}.
The following special case is known by Hesselholt and Madsen
\cite[Thm. A]{larsmad}:

\begin{theorem}\label{dvfield}
Let $K$ be a complete discrete valuation field of
characteristic $0$ with perfect residue field of 
characteristic $p>2$. Then for $i\geq 1$,
$$K_i(K,\Z/p^r)\cong K_i^\et(K,\Z/p^r).$$
\end{theorem}

This has been generalized to certain discrete valuation rings with
non-perfect residue fields in \cite{goodpaper}.
See \cite{larsicm} for a survey of these results.

\subsection{Topological cyclic homology}
Using the hyper-cohomology construction of Thomason, 
one can \cite{ichlars} extend the definition of 
topological Hochschild homology $\TH(A)$ for
a ring $A$ by considering the presheaf of spectra
$\TH:U\mapsto \TH(\Gamma(U,{\mathcal O}_U))$, and setting
\begin{equation}\label{thdescent}
\TH(X_\tau)={\mathbb H}^\cdot(X_\tau,\TH).
\end{equation}

\begin{proposition}\label{etiszar}
a) \cite[Cor. 3.3.3]{ichlars} If the Grothendieck topology $\tau$ 
on the scheme $X$ is coarser than or equal to
the \'etale topology, then $\TH(X_\tau)$ is independent of
the topology (and we drop $\tau$ from the notation).

b) \cite[Cor. 3.2.2]{ichlars} If $X$ is the spectrum of a ring $A$, then 
$\TH(A)\stackrel{\sim}{\longrightarrow}\TH(X)$.
\end{proposition}

It follows from the proposition and \eqref{cechissheaf}
that for a notherian scheme of finite Krull dimension,
\begin{equation}\label{cechdescent}
\TH(X)\cong \TH(X_\Zar)\cong \check {\mathbb H}^\cdot({\mathcal U},\TH).
\end{equation}
In particular, $\TH(X)$ is determined by the spectra
$\TH(U_i)$ for $U_i\in \{{\mathcal U}\}$.

The spectrum $\TR^m(X;p)$ is the fixed point spectrum under of the 
cyclic subgroup of roots of unity $\mu_{p^{m-1}}\subseteq S^1$
acting on $\TH(X)$; let $\TR^m(X;p,\Z/p^r)$ be the version with
coefficients.
There are natural maps called Frobenius and restriction map
$$F,R:\TR^m(X;p,\Z/p^r)\to \TR^{m-1}(X;p,\Z/p^r),$$
and topological cyclic homology 
$\TC^m(X;p,\Z/p^r)$ is the homotopy equalizer of $F$ and $R$.
We view $\TC^\cdot(X;p,\Z/p^r)$ as a pro-spectrum with $R$ as the 
structure map, and define
$$\TC(X;p,\Z_p)=\holim_{m,r} \TC^m(X;p,\Z/p^r).$$
If $(TC^m/p^r)_i$ is the sheaf associated to the presheaf
$U\mapsto \TC_i^m(U;p,\Z/p^r)$, then \eqref{thomsscont} takes the form
\begin{equation}
\label{tcspectrals}
E_2^{s,t}=H^s_{cont}(X_\et,(TC^\cdot/p^\cdot)_{-t}) \Rightarrow
\TC_{-s-t}(X;p,\Z_p).
\end{equation}
If we use the Zariski or Nisnevich topology instead of the 
\'etale topology, we get a different spectral sequence
with the same abutment. The statements of Proposition \ref{etiszar}
and \eqref{cechdescent} procreate to analog statements for $\TC$.

Topological cyclic homology comes equipped with 
the cyclotomic trace map
$$\tr':K(X,\Z_p)\to \TC(X;p,\Z_p).$$
In \cite{larsmad1}, Hesselholt and Madsen 
show that the trace map is an isomorphism
in non-negative degrees for a finite algebra
over the Witt ring of a perfect field. 
Since Thomason's construction is functorial, this factors by 
Proposition \ref{etiszar} a) through
$$\tr:K^\et(X,\Z_p)\to \TC(X_\et;p,\Z_p)\cong \TC(X;p,\Z_p).$$

\begin{theorem}\label{maintool}
\cite[Thm. A]{ichlarsmixed}
Let $X$ is a smooth, proper scheme over a henselian
discrete valuation ring $V$ of mixed characteristic $(0,p)$.
Then the cyclotomic trace
map from \'etale $K$-theory to topological cyclic homology 
\begin{equation*}
K_i^\et(X,\Z_p)\stackrel{\tr_i}{\longrightarrow}\TC_i(X;p,\Z_p)
\end{equation*}
is an isomorphism. 
\end{theorem}

\section{The trace map for arithmetic schemes}
We fix a prime $p\not=2$, let $\Z_h$ the henselization, and
$\Z_p$ the completion of the integers at $p$,
$\Q_h=\Z_h[\frac{1}{p}]$, and $\Q_p=\Z_p[\frac{1}{p}]$.
We fix a number field $K$, and let $\O$ be its ring of integers.
For a set of prime ideals $S$ of $\O$ not containing any of the
primes dividing $p$, we let $\O_S$ be the $S$-integers of $K$. 

\begin{proposition}
Let $X$ be a regular scheme, flat and proper over $\Spec \O_S$,
with good reduction
at $p$. Then there is a commutative diagram
$$\begin{CD}
K_i(X,\Z_p) @>\alpha_i>> K_i^\et(X,\Z_p)
@>f_i >> K_i^\et(X\times_\Z\Z_h,\Z_p) \\
@VV \tr_i' V @V \tr_iVV @V\cong VV \\
\TC_i(X;p,\Z_p) @= \TC_i(X;p,\Z_p) @>\sim>>\TC_i(X\times_\Z\Z_h;p,\Z_p) .
\end{CD}$$
In particular, the cyclotomic trace map is isomorphic to the
composition
$$ K_i(X,\Z_p)\stackrel{\alpha}{\longrightarrow} K_i^\et(X,\Z_p)
\stackrel{f_i}{\longrightarrow} K_i^\et(X\times\Z_h,\Z_p).$$ 
The same statements hold with $\Z_p$ instead of $\Z_h$.
\end{proposition}

\begin{proof} Clearly the diagram commutes, and the right vertical map 
is an isomorphism by Theorem \ref{maintool}.
To show the isomorphism
$$\TC_i(X;p,\Z_p)\stackrel{\sim}{\longrightarrow} 
\TC_i(X\times_\Z\Z_h;p,\Z_p),$$
let ${\mathcal U}=\{U_i\}$ be an affine open covering of $X$.
Then ${\mathcal U}_h=\{U_i\times_\Z\Z_h\}$ is an affine open
covering of $X\times_\Z\Z_h$. By \eqref{cechdescent}, it suffices
to show that $\TC(U_i;p,\Z_p)$ is homotopy equivalent to 
$\TC(U_i\times_\Z\Z_h;p,\Z_p)$. Thus we can assume that $X=\Spec R$, 
with $R$ flat and of finite type over $\O_S$. Then $p$ is not a zero
divisor in $R$, and the rings $R/p^s$ and
$(R\otimes\Z_h)/p^s$ are isomorphic. 
By \cite[Addendum 3.1.2]{ichlarsmixed} we get
$ \TC_i(R;p,\Z_p) \stackrel{\sim}{\longrightarrow}
\TC_i(R\otimes\Z_h;p,\Z_p)
\stackrel{\sim}{\longrightarrow} \TC_i(R\otimes \Z_p;p,\Z_p)$.
\end{proof}

The problem in evaluating the
trace map $\tr_i$ is the calculation of the \'etale
$K$-theory groups involved. We solve this problem by 
localizing away from $p$:

\begin{theorem}
\label{uiop}
Let $X$ be regular scheme, flat and proper over $\Spec \O_S$,
with good reduction at $p$.

a) If 
$K_i^\et(X\times_\Z\Z[\frac{1}{p}],\Z_p)\stackrel{g_i}{\longrightarrow}
K_i^\et(X\times_\Z\Q_p,\Z_p)$ is the restriction map, then there is an 
exact sequence
$$0\to \coker \tr_i\to \coker g_i \stackrel{\delta}{\longrightarrow}
\ker \tr_{i-1}\to \ker g_{i-1}\to 0.$$

b) Let $d$ be the relative dimension of $X$ over $\O_S$.
If $i\geq \frac{8}{3}(d+2)(d+3)(d+4)-14$, or if
the Lichtenbaum-Quillen conjecture holds and $i>d+1$, 
then the map $\delta$ is the zero map.
\end{theorem}

\begin{proof} a) For a closed subset $Z$ of a scheme $X$ with open
complement $U$, we let $K^{\et,Z}(X,\Z_p)$ be the homotopy fiber of
the natural map $K^\et(X,\Z_p)\to K^\et(U,\Z_p)$. 
The closed complement $Y=X\times_{\Z} \F_p$ of
$X\times_\Z\Z[\frac{1}{p}]$ in $X$
is isomorphic to the closed complement of
$X\times\Q_h$ in $X\times\Z_h$.
Consider the natural map of long exact sequences
\begin{equation}\label{dag}
\begin{CD}
K_i^{\et,Y}(X,\Z_p)@>>>
K_i^\et(X,\Z_p)@>>> K_i^\et(X\times_\Z\Z[\frac{1}{p}],\Z_p) \\
@VVV @Vf_iVV @Vg_i^hVV \\
K_i^{\et,Y}(X\times\Z_h,\Z_p)@>>> K_i^\et(X\times\Z_h,\Z_p)@>j^*>>
K_i^\et(X\times\Q_h,\Z_p).
\end{CD}
\end{equation}
By Theorem \ref{maintool}, $f_i$ can be identified
with $\tr_i$, and by the following Lemma $g_i$ can be identified
with $g_i^h$. According to \cite[Thm. D.4]{tt}, there are spectral
sequences
$$\begin{CD}
E_2^{s,t}=H^s_Y(X,({\mathcal K}/p^r)_{-t}) @E\Rightarrow EE  
K_{-s-t}^{\et,Y}(X,\Z/p^r)\\
@VVV @VVV \\
E_2^{s,t}=H^s_Y(X\times\Z_h,({\mathcal K}/p^r)_{-t}) @E\Rightarrow EE
K_{-s-t}^{\et,Y}(X\times \Z_h,\Z/p^r).
\end{CD}$$
By \'etale excision \cite[Prop. 1.27]{milneetale}, the
$E_2$-terms of the two spectral sequences are isomorphic, because
$X\times_\Z\Z_h$ is the direct limit of \'etale neighborhoods
of $Y$ in $X$. Taking the
limit over $r$ shows that the two left terms in diagram \eqref{dag} are 
isomorphic, and we get a) by an easy diagram chase. 

b) It suffices to show that the map $j^*$
in diagram \eqref{dag} is surjective.
Consider the commutative diagram
$$\begin{CD}
K_i(X\times\Z_h,\Z_p)@>>> K_i(X\times\Q_h,\Z_p) \\
@VVV @VVV \\
K_i^\et(X\times\Z_h,\Z_p)@>>>K_i^\et(X\times\Q_h,\Z_p).
\end{CD}$$
The right hand map is surjective for 
$i\geq \frac{8}{3}(d+2)(d+3)(d+4)-14$ by Thomason 
\cite{thomasonsurjective}, and
the Lichtenbaum-Quillen conjecture implies that
the right hand map is an isomorphism for $i>d+1$. 
On the other hand, by localization,
the cokernel of the upper map is contained in 
$K_{i-1}(Y,\Z_p)$, which is zero for $i-1>d$
by \cite{marcI} because $Y$ is smooth. 
Hence the lower map is surjective
\end{proof}

\begin{lemma}
\label{hghg}
Let $X$ be a smooth scheme over $\Spec \Q_h$. Then for any $i$, we
have
$$K_i^\et(X,\Z_p) \cong K_i^\et(X\times_{\Q_h} \Q_p,\Z_p).$$
\end{lemma}

\begin{proof} In view of spectral sequence \eqref{qwer},
$$\begin{CD}
\textstyle E_2^{s,t}=H^s(X,\Z_p(-\frac{t}{2})) @E\Rightarrow EE
K_{-s-t}^\et(X,\Z_p) \\
@VVV @VVV \\
\textstyle E_2^{s,t}=H^s(X\times \Q_p,\Z_p(-\frac{t}{2})) 
@E\Rightarrow EE K_{-s-t}^\et(X\times_{\Q_h}{\Q_p},\Z_p),
\end{CD}$$
it suffices to show that the canonical map induces an
isomorphism on $E_2$-terms. By \cite[Cor. 3.4]{jan} there are
spectral sequences of \'etale cohomology groups
$$\begin{CD}
E_2^{s,t}=H^a(\Q_h,H^b(X\times \bar\Q_h,\Z_p(n))) 
@E \Rightarrow EE
H^{a+b}(X,\Z_p(n))\\
@VVV @VVV \\
E_2^{s,t}=H^a(\Q_p,H^b(X\times \bar\Q_p,\Z_p(n))) @E\Rightarrow EE
H^{a+b}(X\times \Q_p,\Z_p(n)).
\end{CD}$$
The Galois groups of $\Q_h$ and $\Q_p$ are isomorphic, and so are
the Galois modules. Indeed, this is a consequence of the smooth base
change theorem for finite coefficients \cite[Cor. VI 4.3]{milneetale}, 
and immediately extends to continuous cohomology. 
\end{proof}

\begin{corollary}
Let $X$ be regular scheme, flat and proper over 
$\Spec \O_S$ with good reduction at the primes above $p$. 
Assume that the Lichtenbaum-Quillen conjecture
holds, that $i>\dim X$ and $p>\dim X+2$. Then the map
$\tr_i: K_i^\et(X,\Z_p) \to \TC_i(X,\Z_p)$ is isomorphic to 
the sum of localization maps
$$\bigoplus_{a}{\textstyle 
H^{2a-i}(X\times_\Z\Z[\frac{1}{p}],\Z_p(a))}\to
\bigoplus_{a} H^{2a-i}(X\times_\Z\Q_p,\Z_p(a)).$$
\end{corollary}

\begin{proof} The spectral sequences
$$\begin{CD}
\textstyle
E_2^{s,t}=H^s(X\times_\Z\Z[\frac{1}{p}],\Z_p(-\frac{t}{2})) 
@E\Rightarrow EE
K_{-s-t}^\et(X\times_\Z\Z[\frac{1}{p}],\Z_p)\\
@VVV @VVV \\
$$\textstyle 
E_2^{s,t}=H^s(X\times_\Z\Q_p,\Z_p(-\frac{t}{2})) @E \Rightarrow EE
K_{-s-t}^\et(X\times_\Z\Q_p,\Z_p)
\end{CD}$$
degenerate at $E_2$ with split filtration for
$p>\frac{cd_p X}{2}$ by Soul\'e \cite[Thm. 1]{souleappl}. 
\end{proof}

The localization map for \'etale cohomology is related to
$L$-functions and $p$-adic $L$-functions by Iwasawa theory. In the
following sections, we give examples for number fields, in
particular totally real number fields. It should be possible to
extend these results to Dirichlet characters, elliptic curves with
complex multiplication, or Hecke characters of imaginary quadratic
fields as in \cite{soulepadic, thesis}.

\section{Number fields}
In the case of a number field, we can make the calculations of the
last section more explicit. The translation of results of Iwasawa-theory
into results on regulators are similar to \cite{kolster}.
We keep the notation of the previous section.
Let $\p$ be a prime of $\O_S$ dividing $p\not=2$,
$\O_\p$ be the completion of $\O_S$ at $\p$,
$K_\p$ its quotient field, and $k_\p\cong \O_\p/{\mathfrak m_\p}$
its residue field. Similarly, let $\O_\p^h$ be the henselization of 
$\O_S$ at $\p$, and $K_\p^h$ its quotient field. The residue fields
of $\O_\p^h$ and $\O_\p$ are canonicallly isomorphic. 
Note that 
$K\otimes_\Z\Q_p\cong\prod_{{\mathfrak p}|p}K_{\mathfrak p}$, 
$\O_S\otimes_\Z\Z_p\cong\prod_{{\mathfrak p}|p}\O_{\mathfrak p}$, and 
$(\O_S\otimes_\Z\F_p)^{red}\cong \prod_{{\mathfrak p}|p}k_\p$,
and similarly for the henselization.

\begin{proposition}\label{plpl}
For $i>1$ and $j=1,2$ we have the following isomorphisms
\begin{align*}
K_{2i-j}^\et(\O_S,\Z_p)\stackrel{\sim}{\longrightarrow}
&\textstyle K_{2i-j}^\et(\O_S[\frac{1}{p}],\Z_p)
\stackrel{\sim}{\longrightarrow}
\textstyle H^j(\Spec \O_S[\frac{1}{p}],\Z_p(i))\\
K_{2i-j}^\et(\O_\p,\Z_p)\stackrel{\sim}{\longrightarrow}
&K_{2i-j}^\et(K_\p,\Z_p)
\stackrel{\sim}{\longrightarrow} H^j(\Spec K_\p, \Z_p(i)).
\end{align*}
\end{proposition}

\begin{proof} 
Since $\Spec \O_S[\frac{1}{p}]$ and $K_\p$ have cohomological dimension 
$2$ if $p\not=2$, the right hand
isomorphism follows from the spectral sequence \eqref{qwer} and 
$H^0(\Spec \O_S[\frac{1}{p}],\Z_p(i))=H^0(\Spec K_\p, \Z_p(i))=0$
for $i>0$.

Consider the commutative diagram of long exact sequences
$$\begin{CD}
K_i(k_\p,\Z_p)@>>> K_i(\O_\p,\Z_p)@>>> K_i(K_\p,\Z_p)\\
@VVV @VVV @VVV \\
K_i^{\et,k_p}(\O_\p,\Z_p)@>>>
K_i^\et(\O_\p,\Z_p)@>>> K_i^\et(K_\p,\Z_p)\\
\end{CD}$$
The middle map is an isomorphism by \cite{larsmad1} and the
the right map is an isomorphism by Theorem \ref{dvfield}, 
hence in view of $K_i(k_\p,\Z_p)=0$ for $i>0$ we get
the local result together with $K_i^{\et,k_p}(\O_\p,\Z_p)=0$ for $i>0$.
Comparing the lower row with the analog row for
the henselization, we get $K_i^{\et,k_\p}(\O_\p^h,\Z_p)=0$ for $i>0$.
Indeed, $K_i^\et(K_\p^h,\Z_p)\cong K_i^\et(K_\p,\Z_p)$ by
Lemma \ref{hghg}, and by \eqref{contin} there are isomorphisms 
$$\begin{CD}
K_i(\O_\p^h,\Z_p)@>\sim >> K_i(\O_\p,\Z_p)\\
@V\cong VV @V\cong VV\\
K_i^\et(\O_\p^h,\Z_p)@> >> K_i^\et(\O_\p,\Z_p),
\end{CD}$$
Using \'etale excission as in the proof of Theorem \ref{uiop} a), we see
that 
$K_i^{\et,k_\p}(\O_S,\Z_p)\cong K_i^{\et,k_\p}(\O_\p^h,\Z_p)=0$, 
hence the global result.
\end{proof}


Let $\Sigma$ be the set of primes of $\O$ dividing $p$ or
infinity, and
$G_\Sigma$ be the Galois group of the maximal extension
of $K$ unramified outside $\Sigma$. Then by 
\cite[\S 3.2]{coateslicht},
$$\textstyle
H^n(\Spec \O[\frac{1}{p}],\Z_p(i))\cong 
H^n_{cont}(G_\Sigma, \Z_p(i)),$$
where the left hand side is \'etale cohomology and the right
hand side is continuous Galois cohomology.
Proposition \ref{plpl} for $S=\emptyset$ 
shows that for $i>0$ the trace map can be identified:
$$\begin{CD}
K_{2i-j}^\et(\O,\Z_p) @>tr_{2i-j}>> \TC_{2i-j}(\O;p,\Z_p) \\
@| @| \\
H^j(G_\Sigma,\Z_p(i)) @>>> \prod_{\p|p} H^j(K_\p,\Z_p(i)).
\end{CD} $$

The following fundamental conjecture is due to Schneider 
\cite[p. 129]{schneider}:
\medskip

\noindent {\bf Conjecture S(K,i)}
\label{schn}  {\it The group 
$ H^2(G_\Sigma,\Z_p(i))$ is torsion.}{\bf .}\medskip

By Soul\'e \cite{soule}, $S(K,i)$ holds for all $i>1$, and 
$S(K,i)$ is by the long exact coefficient sequence equivalent to 
$H^2(G_\Sigma,\Q_p/\Z_p(i))=0$.

\begin{lemma}
\label{cohomology} Let $K$ be a number field of degree $d$ over
$\Q$ with $r_1$ real and $r_2$ complex embeddings.

a) \cite[Satz 3.2, 3.4]{schneider} If $i\not=0,1$, then
$ \rk_{\Z_p} \prod_{\p|p}H^1(K_\p,\Z_p(i)) =d$, and the group 
$H^2(K_\p,\Z_p(i))$ is finite.

b) \cite[Satz 4.6]{schneider}
$$ \rk H^1(G_\Sigma,\Z_p(i)) =\begin{cases}
r_2+ \rk H^2(G_\Sigma,\Z_p(i)) & i\not= 0 \text{ even};\\
r_1+r_2 +\rk  H^2(G_\Sigma,\Z_p(i))& i\not= 1 \text{ odd}.
\end{cases}$$
\end{lemma}

For $i>0$, the Tate-Poitou exact sequence \cite{milneadt}
\begin{multline}
\label{tp}
0 \to H^2(G_\Sigma,\Q_p/\Z_p(1-i))^* \to \\
H^1(G_\Sigma,\Z_p(i)) \stackrel{\tr_{2i-1}}{\longrightarrow}
\prod_{\p|p}H^1(K_\p,\Z_p(i))
\to H^1(G_\Sigma,\Q_p/\Z_p(1-i))^* \to \\
H^2(G_\Sigma,\Z_p(i)) \stackrel{\tr_{2i-2}}{\longrightarrow}
\prod_{\p|p}H^2(K_\p,\Z_p(i)) \to 
H^0(G_\Sigma,\Q_p/\Z_p(1-i))^* \to 0.
\end{multline}
can be used to get results on the trace map.

For example, $S(K,1-i)$ is equivalent to
the injectivity of $tr_{2i-1}$, and the kernels of cokernels
of the trace maps can be expressed in terms of cohomology groups:
\begin{align}\label{corir}
\ker \tr_{2i-1}&\cong H^2(G_\Sigma,\Q_p/\Z_p(1-i))^*\\
0\to \coker \tr_{2i-1} \to &H^1(G_\Sigma,\Q_p/\Z_p(1-i))^*
\to \ker \tr_{2i-2}\to 0\\
\coker \tr_{2i-2} &\cong H^0(G_\Sigma,\Q_p/\Z_p(1-i))^*
\end{align}
Moreover, by \cite[Satz 5 ii]{schneider}
\begin{multline*}
\ker tr_{2i-2}=
\coker  \big( H^1(K_p,\Z_p(i)) \to H^1(G_\Sigma,\Q_p/\Z_p(1-i))^*\big) \\
=\ker \big(H^1(G_\Sigma,\Q_p/\Z_p(1-i))\to 
\prod_{\p|p}H^1(K_\p,\Q_p/\Z_p(1-i))^*\big) \\
=\frac{\div H^1(K,\Q_p/\Z_p(i))}{\Div H^1(K,\Q_p/\Z_p(i))} .
\end{multline*}

\subsection{Totally real fields}
Let $K$ be a totally real field of degree $d$ over $\Q$. Let $w_i$ and
$w_{{\mathfrak p},i }$ be the order of the group
$H^0(G_\Sigma,\Q_p/\Z_p(j))$ and $H^0(K_{\mathfrak p},\Q_p/\Z_p(j))$,
respectively. We normalize the $p$-adic absolute value such that
$|p|_p=\frac{1}{p}$, and write equalities of $p$-powers so
that only positive powers appear on both sides of an equation.

\begin{lemma}
\label{realcohom} a) Let $i>0$ be even. Then the groups
$H^0(G_\Sigma,\Q_p/\Z_p(i))\cong H^1(G_\Sigma,\Z_p(i))$,
and $H^1(G_\Sigma,\Q_p/\Z_p(i)) \cong H^2(G_\Sigma,\Z_p(i))$
are finite.

b) Let $i\not=1$ be odd. Then $H^0(G_\Sigma,\Q_p/\Z_p(i))=0$,
$H^1(G_\Sigma,\Z_p(i))\cong \Z_p^d$, and $S(K,i)$ holds.

c) If $K/\Q$ is unramified at $p$ and $p-1 \not| i$, then for
every ${\mathfrak p}|p$,
$$H^0(K_{\mathfrak p},\Q_p/\Z_p(i))= H^2(K_{\mathfrak p},\Z_p(1-i))=0.$$
\end{lemma}

\begin{proof} a) By $S(K,i)$ for $i>1$ and Lemma \ref{cohomology} b), 
we get that $H^1(G_\Sigma,\Q_p(i))=H^2(G_\Sigma,\Q_p(i)=0$ for $i>0$ even.
The result now follows from the long exact coefficient sequence.

b) Since $K$ is totally real and $p\not=2$, the
degree of the extension
$K(\mu_p)/K$ is divisible by two, hence
$|H^0(G_\Sigma,\Q_p/\Z_p(i))|=\max\{p^j:[K(\mu_{p^j}):K]|i\}=1$.
As this group is isomorphic to the torsion subgroup of $H^1(G_\Sigma,\Z_p(i))$,
the latter group is torsion free. The group $H^2(G_\Sigma,\Z_p(i))$
is finite because $H^2(G_\Sigma,\Q_p/\Z_p(i))$ 
is zero. For $i>1$ this follows from $S(K,i)$, and for $i<0$ its dual is a
subgroup of $H^1(G_\Sigma,\Z_p(1-i))$ by \eqref{tp}, and the latter
is finite by a). 

c) Since the extension $K_{\mathfrak p}/\Q_p$ is unramified by hypothesis.
But $\Q_p(\mu_p)/\Q_p$ is totally ramified at $p$ and has degree $p-1$, 
hence the same holds for the
extension $K_{\mathfrak p}(\mu_p)/K_{\mathfrak p}$. 
We can now use the argument of b) together with local duality
\cite[Satz 2.4]{schneider}
\end{proof}

\begin{proposition}
Let $i>0$ be even. Then
$\ker \tr_{2i-1}=\coker \tr_{2i-2} = 0$, and
$$  |\zeta_K(1-i)|_p^{-1}\cdot|\ker \tr_{2i-2}| \cdot 
\prod_{\mathfrak p}w_{{\mathfrak p},1-i}= w_{i}.$$
\end{proposition}

\begin{proof} The first two statements follow from Lemma \ref{realcohom} b),
and \eqref{corir}. For the zeta value, we have
\begin{multline*}
 |\zeta_K(1-i)|_p =
\frac{|H^0(G_\Sigma,\Q_p/\Z_p(i))|}
{|H^1(G_\Sigma,\Q_p/\Z_p(i))|}
=\frac{|H^1(G_\Sigma,\Z_p(i))|}
{|H^2(G_\Sigma,\Z_p(i))|} \\
=\frac{|\ker \tr_{2i-2}|\cdot \prod_{\p|p}|H^2(K_\p,\Z_p(i))|}{w_i} =
\frac{|\ker \tr_{2i-2}| \cdot \prod_{\mathfrak p}w_{{\mathfrak p},1-i}}{w_{i}}.
\end{multline*}
The first equality is \cite[Thm. 6.2]{bn}, the second equality
is Lemma \ref{realcohom} a), the third equality
follows from \eqref{tp} because $\tr_{2i-2}$ is surjective, and the 
last equality is local duality. 
\end{proof}

\begin{theorem}
Let $i>0$ is odd. Then $S(K,1-i)$ holds if and only if $\tr_{2i-1}=0$
if and only if $\tr_{2i-1}$ has finite cokernel. In this case,
$$|\coker \tr_{2i-1}|\cdot|H^2(G_\Sigma,\Z_p(i))|
=\prod_{{\mathfrak p}|p}w_{{\mathfrak p},1-i}\cdot
|L_p(K,\omega^{1-i},i)|_p^{-1}.$$
\end{theorem}

\begin{proof} If $i>1$ is odd, then by Lemma \ref{realcohom} b)
$$\rk H^1(G_\Sigma, \Z_p(i)) = \rk \prod_{{\mathfrak p}|p}
 H^1(K_{\mathfrak p}, \Z_p(i)) =d.$$
By the Tate-Poitou sequence \eqref{tp},
$H^2(G_\Sigma,\Q_p/\Z_p(1-i))=0$ if and only if
$\ker \tr_{2i-1}=0$ if and only if
$\tr_{2i-1}$ has finite cokernel. In this case,
by \cite[Thm. 6.1]{bn}, and \eqref{tp},
\begin{multline*}
|L_p(K,\omega^{1-i},i)|_p= \frac{| H^0(G_\Sigma,\Q_p/\Z_p(1-i))|}
{|H^1(G_\Sigma,\Q_p/\Z_p(1-i))|} \\
=\frac{\prod_{{\mathfrak p}|p} |H^2(K_{\mathfrak p},\Z_p(i))| }
{|\coker \tr_{2i-1}| \cdot |H^2(G_\Sigma,\Z_p(i))|} .
\end{multline*}
The result follows with local duality.
\end{proof}

The first result in this direction is due to 
Soul\'e \cite{soulepadicreg}. He maps a subgroup of the source
of $\tr_{2i-1}$ to a quotient of the target, and relates the
index of this map to the $p$-adic $L$-function directly 
(without using the main theorem of Iwasawa theory).

Note the difference between the case $i$ even and $i$ odd: In the
former case, the Euler characteristic of $\Q_p/\Z_p(i)$ gives a a
result on the $p$-adic $L$-function at $1-i$ which translates into a
result for the $\zeta$-function at $1-i$, because the $p$-adic
$L$-function approximates the $\zeta$-function at negative integers. In
the latter case, the Euler characteristic of
$\Q_p/\Z_p(1-i)$ only gives a result for the $p$-adic $L$-functions
at $i$. 

We give a version for $K$-theory instead of \'etale
$K$-theory:

\begin{corollary} Assume $i>1$ is an odd integer.

a) The trace map factors like
$$\tr_{2i-1}'  :K_{2i-1}(\O)\otimes \Z_p/\text{tors} \to 
\TC_{2i-1}(\O,\Z_p) .$$

b) $S(K,1-i)$ implies
$$ |\coker \tr_{2i-1}'|\cdot |H^2(G_\Sigma,\Z_p(i))|=
\prod_{{\mathfrak p}|p}w_{{\mathfrak p},1-i}\cdot
|L_p(K,\omega^{1-i},i)|_p^{-1}.$$

c) If moreover $K/\Q$ is unramified at $p$ and $i\not\equiv 1\mod p-1$,
then
$\tr_{2i-1}'$ is a map between free $\Z_p$-modules of rank
$d$, and
$$|\coker tr_{2i-1}'|\cdot |H^2(G_\Sigma,\Z_p(i))|=
|L_p(K,\omega^{1-i},i)|_p^{-1}.$$

d) Let $K=\Q$, $i\not\equiv 1\mod p-1$, and
$n$ a positive integer with $n\equiv -i \mod p-1$.
If $p| {\frac{B_{n+1}}{n+1}}$, then
$H^2(G_\Sigma,\Z_p(i))\not=0$ or $\tr'_{2i-1}$ is not surjective.
\end{corollary}

\begin{proof} 
a) By Lemma \ref{realcohom} b) and Proposition \ref{plpl}
$$K_{2i-1}^\et (\O,\Z_p)_{tors} =
H^1(G_\Sigma,\Z_p(i))_{tors}= H^0(G_\Sigma, \Q_p/\Z_p(i))=0,$$
so the trace map factors through the torsion free quotient
of $K_{2i-1}(\O)\otimes \Z_p$.

b) By Soul\'e \cite{soule} the map
$K_{2i-1}(\O,\Z_p)/{tors}\to K_{2i-1}^\et(\O,\Z_p)/{tors}$
is surjective. Since both groups are free $\Z_p$-modules of the
same rank, they must be isomorphic. In other words,
$K_{2i-1}(\O,\Z_p)/tors \cong H^1(G_\Sigma,\Z_p(i))$, and the
statement is a just a reformulation of the previous theorem.

c) The torsion subgroup of $\TC_{2i-1}(\O;p,\Z_p)\cong 
\prod_{\p|p}H^1(K_\p,\Z_p(i))$
is isomorphic to $\prod_{\p|p}H^0(K_\p,\Q_p/\Z_p(i))=0$ 
by Lemma \ref{realcohom} c).

d) By \cite[Thm. 5.11, Cor. 5.13]{washington},
\begin{multline*}
L_p(\Q,\omega^{1-i},i)=L_p(\Q,\omega^{n+1},i)\equiv 
L_p(\Q,\omega^{n+1},1-(n+1))\\
\equiv -(1-p^n)\frac{B_{n+1}}{n+1}\equiv -\frac{B_{n+1}}{n+1} \mod
p .
\end{multline*}
\end{proof}


\begin{thebibliography}{99}

\bibitem{bn} {\sc P.Bayer, J.Neukirch}, On values of zeta functions
and $l$-adic Euler characteristics. Inv. Math. {\bf 50} (1978),
35--64.












\bibitem{browngersten} {\sc K.S.Brown, S.M.Gersten}, 
Algebraic $K$-theory as generalized sheaf cohomology. 
Algebraic K-theory, I: Higher K-theories (Proc. Conf., Battelle Memorial Inst.
Seattle, Wash., 1972), pp. 266--292. Lecture Notes in Math., Vol. {\bf 341}, 
Springer, Berlin, 1973. 

\bibitem{coateslicht} {\sc J.Coates, S.Lichtenbaum}, 
On $l$-adic zeta functions. Ann. of Math. (2) {\bf 98} (1973), 498--550.








\bibitem{gabberrig} {\sc O.Gabber}, $K$-theory of Henselian local
rings and Henselian pairs. Algebraic $K$-theory, commutative
algebra, and algebraic geometry (Santa Margherita Ligure, 1989).
Cont. Math. {\bf 126} (1992), 59--70.


\bibitem{thesis} {\sc T.Geisser}, $p$-adic $K$-theory of Hecke
characters of imaginary quadratic fields. Duke Math. J. {\bf 86}
(1997), 197--238.

\bibitem{handbook} {\sc T.Geisser}, Motivic cohomology, K-theory
and topological cyclic homology. To appear in: K-theory handbook.





\bibitem{ichlars} {\sc T.Geisser, L.Hesselholt}, Topological Cyclic
Homology of Schemes. K-theory, Proc. Symp. Pure Math. AMS {\bf 67}
(1999), 41--87.

\bibitem{ichlarsmixed} {\sc T.Geisser, L.Hesselholt}, $K$-theory and
topological cyclic homology of smooth schemes over discrete
valuation rings. Trans. Amer. Math. Soc. (to appear)

\bibitem{goodpaper} {\sc T.Geisser, L.Hesselholt}, The de Rham-Witt
complex and $p$-adic vanishing cycles. Preprint 2003.

\bibitem{marcI} {\sc T.Geisser, M.Levine}, The $p$-part of
$K$-theory of fields in characteristic $p$. Inv. Math. {\bf 139}
(2000), 459--494.










\bibitem{larsicm} {\sc L.Hesselholt}, Algebraic $K$-theory and trace invariants.
Proceedings of the International Congress of Mathematicians, Vol. II 
(Beijing, 2002), 415--425, Higher Ed. Press, Beijing, 2002.

\bibitem{larsmad1} {\sc L.Hesselholt, I.Madsen}, On the $K$-theory of finite 
algebras over Witt vectors of perfect fields. 
Topology {\bf 36} (1997), no. 1, 29--101.


\bibitem{larsmad} {\sc L.Hesselholt, I.Madsen}, The $K$-theory
of local fields. To appear in: Annals of Math.


\bibitem{jan} {\sc U.Jannsen}, Continuous \'etale cohomology. Math. Ann. {\bf
280} (1988), 207--245.







\bibitem{kolster} {\sc M.Kolster, Nguyen Quang Do, Thong}, Syntomic regulators 
and special values of $p$-adic $L$-functions. 
Invent. Math. {\bf 133} (1998), no. 2, 417--447.











\bibitem{mccarthy} {\sc R.McCarthy}, Relative algebraic K-theory
and topological cyclic homology. Acta Math. {\bf179} (1997), 197--222.

\bibitem{milneadt} {\sc J.S.Milne}, Arithmetic duality theorems. Perspectives 
in Mathematics, 1. Academic Press, Inc., Boston, MA, 1986.



\bibitem{milneetale} {\sc J.S.Milne}, Etale cohomology. Princeton Math.
Series 33.




\bibitem{nisnevich} {\sc Y.Nisnevich}, The completely decomposed topology on
schemes and associated descent spectral sequences in algebraic
$K$-theory. Algebraic $K$-theory: connections with geometry and
topology (Lake Louise, AB, 1987), NATO Adv. Sci. Inst. Ser. C
Math. Phys. Sci. {\bf 279} (1989), 241--342.


\bibitem{panin} {\sc I.A.Panin}, The Hurewicz theorem and $K$-theory of 
complete discrete valuation rings. Izv. Akad. Nauk SSSR Ser. Mat. 
{\bf 50} (1986), no. 4, 763--775.






\bibitem{schneider} {\sc P.Schneider},
\"Uber gewisse Galoiskohomologiegruppen. Math. Zeit. {\bf 168}
(1979), 181--205.


\bibitem{soule} {\sc C.Soul\'e}, $K$-th\'eorie des anneaux d'entiers
de corps de nombres et cohomologie \'etale. Inv. Math. {\bf 55}
(1979), 251--295.

\bibitem{soulepadicreg} {\sc C.Soul\'e}, On higher $p$-adic regulators. 
Algebraic $K$-theory, Evanston 1980 
(Proc. Conf., Northwestern Univ., Evanston, Ill., 1980), pp. 372--401, 
Lecture Notes in Math., {\bf 854}, Springer, Berlin-New York, 1981.


\bibitem{souleappl} {\sc C.Soul\'e}, Operations on etale $K$-theory. 
Applications. Algebraic $K$-theory, Part I (Oberwolfach, 1980), pp. 271--303,
Lecture Notes in Math., {\bf 966}, Springer, Berlin, 1982. 


\bibitem{soulepadic} {\sc C.Soul\'e}, $p$-adic $K$-theory of
elliptic curves. Duke Math. J. {\bf 54} (1987), 249--269.


\bibitem{susdvr} {\sc A.Suslin}, On the K-theory of local fields.
J. Pure Appl. Alg. {\bf 34} (1984), 301--318.




\bibitem{thomasonsurjective} {\sc R.W.Thomason}, Bott stability in algebraic 
$K$-theory. Applications of algebraic $K$-theory to algebraic geometry and 
number theory, Part I, II (Boulder, Colo., 1983), 389--406, Contemp. Math., 
{\bf 55}, Amer. Math. Soc., Providence, RI, 1986.

\bibitem{thomasonetale} {\sc R.W.Thomason}, Algebraic K-theory
and etale cohomology. Ann. Sci. ENS {\bf 18} (1985), 437--552.

\bibitem{tt} {\sc R.Thomason, T.Trobaugh}, Higher algebraic
$K$-theory of schemes and of derived categories.
The Grothendieck Festschrift, Vol. III, 247--435, Progr. Math.,
{\bf 88}, Birkh\"auser Boston, Boston, MA, 1990.




\bibitem{voebk} {\sc V.Voevodsky}, On motivic cohomology with
$\Z/l$-coefficients. Preprint 2003.



\bibitem{washington} {\sc L.Washington}, Introduction to
cyclotomic fields. Springer GTM 83 (1982).


\end{thebibliography}
\end{document}